\documentclass[preprint,11pt]{elsarticle}
\usepackage{amssymb,amsmath,amsthm,amsfonts,amscd}
\usepackage{graphicx,color}
\usepackage[normalem]{ulem} 
\usepackage{soul} 
\usepackage{cancel} 
\usepackage{hyperref}
\definecolor{ForestGreen}{rgb}{0.15,0.416,0.18}
\definecolor{EgyptBlue}{rgb}{0.063,0.2,0.65}
\hypersetup{
	colorlinks=true,
	linkcolor=EgyptBlue,
   citecolor=EgyptBlue,
   urlcolor=ForestGreen
}
\journal{xxxxxxx}
\newtheorem{theorem}{Theorem}[section]

\newtheorem{remark}{Remark}[section]

\newcommand{\cqdf}{\hfill\rule{5pt}{5pt}}

\newcommand{\disp}{\displaystyle}

\def\Om{\Omega}

\def\re{\mathbb{R}}

\def\jnt{\displaystyle\int}
\def\jjnt{\jnt\!\!\!\jnt}

\def\beq{\begin{eqnarray}}
\def\eeq{\end{eqnarray}}
\def\beqa{\begin{eqnarray*}}
\def\eeqa{\end{eqnarray*}}
\def\beqn{\begin{equation}}
\def\eeqn{\end{equation}}
\def\bar{\begin{array}}
\def\ear{\end{array}}
\begin{document}

\begin{frontmatter}

\title{Controllability for a one-dimensional wave equation in a non-cylindrical domain}
\author[IJesus]{Isa\'{i}as Pereira de Jesus\corref{cor1}}
\ead{isaias@ufpi.edu.br}

\address[IJesus]{ Departamento de Matem\'{a}tica, Universidade Federal do Piau\'{i}, 64049-550, Teresina, PI, Brazil}
\cortext[cor1]{Corresponding author. Phone: +55 (86) 3215-5835}

\begin{abstract}
This paper deals with the controllability for a one-dimensional wave equation with mixed boundary conditions in a non-cylindrical domain. This equation models small vibrations of a string where an endpoint is fixed and the other is moving. As usual, we consider one main control (the leader) and an additional secondary control (the follower). We use Stackelberg-Nash strategies.

\end{abstract}


\end{frontmatter}

\section{Introduction and main result}
There are plenty of situations where several controls are required in order to drive a system to one or more objectives. Usually, if we assign different
roles to the controls, we speak of {\it hierarchic control.} The concept of hierarchic control in the context of hyperbolic PDEs was introduced in \cite{L1} when the author analyzed the approximate controllability for a system associated with a wave equation. There, he considered one main control (the leader) and an additional secondary control (the follower).

In \cite{DL,GO}, the hierarchic control of a parabolic PDE and the Stokes systems have been analyzed and used to solve an approximate controllability problem. In \cite{AR, AR300}, a strategy has been used to deduce the exact controllability (to trajectories) for a parabolic PDE. Recently, in \cite {AR501},  the authors studied a hierarchic control problem for a hyperbolic PDE.

In this article, motivated by the arguments contained in the work of J.-L. Lions
\cite{L1}, we investigate a  similar question of hierarchic control for the wave equation, employing the Stackelberg strategy in the case of time
 dependent domains.

In spite of the vast literature on the controllability problems of the wave equation in cylindrical domains, there are only
a few works dealing with non-cylindrical case. We refer to  \cite {Ar, Cui10, Cui, Mi1} for some known results in this direction.

Some other basic references on controllability can be found in \cite {Cor} and \cite {Tucs}.

The novelty of  this paper relies on the consideration of moving boundaries. The approach proposed consists in a suitable change transforming a system written over a moving domain into an equivalent system written over a fixed domain.

\subsection {Statement of the problem}

As in \cite{Cui1}, given $ T > 0 $, we consider the non-cylindrical domain defined
by
$$\displaystyle \widehat{Q}= \left\{ (x,t) \in \re^2;\; 0 < x <  \alpha _k(t);\  t \in (0,T) \right\},$$
where
$$ \alpha_k(t) = 1 + kt, \;\;\;\;\;\; \; 0 < k < 1.$$

Its lateral boundary  is defined by $\displaystyle \widehat{\Sigma} =\widehat{\Sigma}_0
\cup \widehat{\Sigma}_0^*$, with
$$\displaystyle \widehat{\Sigma}_0^* = \{(\alpha_k(t),t);\  t \in (0,T) \} \;\;\;\;\mbox{ and } \;\;\;\; \widehat{\Sigma}_{0} = \widehat{\Sigma} \backslash \widehat{\Sigma}_0^* =  \{(0, t);\  t \in (0,T) \} .$$

We also represent by $\Om_t$ and $\Om_0$ the intervals $\displaystyle (0, \alpha_k(t))$
and $\displaystyle (0, 1)$, respectively.

Consider the following  wave equation in the non-cylindrical domain
$\widehat{Q}$:
\begin{equation} \label{eq1.3}
\left|
\begin{array}{l}
\displaystyle u'' - u_{xx} = 0 \ \ \mbox{ in } \ \ \widehat{Q},\\[11pt]
\displaystyle u(x,t) = \left\{
\begin{array}{l}
\widetilde{w}(t) \ \ \mbox{ on } \ \ \widehat{\Sigma}_{0}^*,\\[11pt]
0 \ \ \mbox{ on } \ \ \widehat{\Sigma}_0,
\end{array}
\right.\\[11pt]
\displaystyle u(x,0) = u_0(x), \;\; u'(x,0) = u_1(x) \ \mbox{ in }\; \Om_0,
\end{array}
\right.
\end{equation}
where $u$ is the state variable, $\widetilde{w}$ is the control variable and $( u_0(x),
u_1(x)) \in L^2(0,1) \times H^{-1}(0,1)$. By $\displaystyle u'=u'(x,t)$ we represent the
derivative $\displaystyle \frac{\partial u}{\partial t}$ and by  $\displaystyle u_{xx}=u_{xx}(x,t)$ the second order partial derivative $\displaystyle \frac{\partial^2 u}{\partial x^2}$. Equation \eqref{eq1.3} models the motion of a string where an endpoint is fixed and the other one is moving. The constant $k$ is called the speed of the moving endpoint.

The main goal of this article is to analyze the hierarchic control of \eqref{eq1.3} and, in particular, to prove that the Stackelberg-Nash strategy allows to solve the approximate controllability problem.

The control system of this paper is similar to that of \cite {Je1000}. But motivated by \cite {Cui1}, the control is put on a different boundary. To overcome the difficulties, we transform the system \eqref{eq1.3} into an equivalent wave equation with variable coefficients in the cylindrical domain and establish the approximate controllability of this equation by Holmgren's Uniqueness Theorem.



Following the work of J.-L. Lions \cite{L1}, we divide $\displaystyle \widehat{\Sigma}_0^*$
into two parts
\begin{equation}\label{decomp0}
\displaystyle \widehat{\Sigma}_0^*=\displaystyle \widehat{\Sigma}_1 \cup \displaystyle \widehat{\Sigma}_2,
\end{equation}
 and consider
\begin{equation} \label{decomp}
\displaystyle \widetilde{w}=\{\widetilde{w}_1, \widetilde{w}_2\}, \;\; \widetilde{w}_i=\mbox{control function in } \; L^2(\widehat{\Sigma}_i), \;i=1,2.
\end{equation}

Thus, we observe  that the system \eqref{eq1.3} can be rewritten as follows:
\begin{equation} \label{eq1.3.1}
\left|
\begin{array}{l}
\displaystyle u''- u_{xx} = 0 \ \ \mbox{ in } \ \ \widehat{Q},\\[5pt]
\displaystyle u(x,t) = \left\{
\begin{array}{l}
\widetilde{w}_1(t) \ \ \mbox{ on } \ \ \widehat{\Sigma}_{1},\\
\widetilde{w}_2(t) \ \ \mbox{ on } \ \ \widehat{\Sigma}_{2},\\
0 \ \ \mbox{ on } \ \ \widehat{\Sigma}\backslash \widehat{\Sigma}_0^*,
\end{array}
\right.\\[5pt]
\displaystyle u(x,0) = u_0(x), \;\; u'(x,0) = u_1(x) \ \mbox{ in }\; \Om_0.
\end{array}
\right.
\end{equation}
In the decomposition \eqref{decomp0}, \eqref{decomp} we establish a hierarchy.
We think of $\widetilde{w}_1$ as being the ``main" control, the leader,  and we
think of $\widetilde{w}_2$ as the follower,  in Stackelberg terminology.

Let us define the following secondary cost functional
\begin{equation}\label{sfn}
\displaystyle \widetilde{J}_2(\widetilde{w}_1, \widetilde{w}_2) =
\displaystyle\frac{1}{2} \displaystyle \jjnt_{\widehat{Q}} \left(u(\widetilde{w}_1, \widetilde{w}_2)-\widetilde{u}_2\right)^2 dx dt +
\displaystyle\frac{\widetilde{\sigma}}{2} \int_{\widehat{\Sigma}_2} \widetilde{w}_2^2\;d\widehat{\Sigma},
\end{equation}
and the main cost functional
\begin{equation}\label{mfn}
\displaystyle \widetilde{J}(\widetilde{w}_1) =\frac{1}{2}\int_{\widehat{\Sigma}_1} \widetilde{w}_1^2\;d\widehat{\Sigma},
\end{equation}
where $\widetilde{\sigma}>0$ is a constant and $\widetilde{u}_2$ is a given
function in $L^2(\widehat{Q}).$

\begin{remark}\label{bdfnc} From the regularity and uniqueness of the solution  of \eqref{eq1.3.1} (see Remark \eqref{rsol}) the cost functionals $\displaystyle \widetilde{J}_2$ and $\displaystyle \widetilde{J}$ are well defined.
\end{remark}

Now, let us describe the Stackelberg-Nash strategy. Thus, for each choice of the leader  $ \widetilde{w}_1,$ we try to find a Nash equilibrium for the cost functional $\widetilde{J}_{2}$, that is, we look for a control $\displaystyle \widetilde{w}_2=\mathfrak{F}(\widetilde{w}_1),$
depending on $\displaystyle \widetilde{w}_1$, satisfying:
\begin{equation}\label{son}
\displaystyle \widetilde{J}_2(\widetilde{w}_1, \widetilde{w}_2)= \inf_{\widehat{w}_2 \in L^2(\widehat{\Sigma}_2)} \widetilde{J}_2(\widetilde{w}_1, \widehat{w}_2).
\end{equation}

After this, we consider the state $\displaystyle u\left(\widetilde{w}_1,
\mathfrak{F}(\widetilde{w}_1)\right)$ given by the solution of
\begin{equation} \label{eq1.3.1F}
\left|
\begin{array}{l}
\displaystyle u''- u_{xx} = 0 \ \ \mbox{ in } \ \ \widehat{Q},\\[3pt]
\displaystyle u(x,t) = \left\{
\begin{array}{l}
\widetilde{w}_{1} \ \ \mbox{ on } \ \ \widehat{\Sigma}_{1},\\
\mathfrak{F}(\widetilde{w}_1) \ \ \mbox{ on } \ \ \widehat{\Sigma}_{2},\\
0 \ \ \mbox{ on } \ \ \widehat{\Sigma}\backslash \widehat{\Sigma}_0^*,
\end{array}
\right.\\[3pt]
\displaystyle u(x,0) = u_0(x), \;\; u'(x,0) = u_1(x) \ \mbox{ in }\; \Om_0.
\end{array}
\right.
\end{equation}

We will look for any optimal control $\displaystyle \widetilde{w_1}$ such that
\begin{equation}\label{ocn}
 \displaystyle \widetilde{J}(\widetilde{w}_{1}, \mathfrak{F}(\widetilde{w}_1))= \inf_{\overline{w}_1 \in L^2(\widehat{\Sigma}_1)} \widetilde{J}(\overline{w}_1, \mathfrak{F}(\widetilde{w}_1)),
\end{equation}
subject to the following restriction of the approximate controllability type
\begin{equation}\label{apcn}
\displaystyle \left(u(x, T; \widetilde{w}_{1}, \mathfrak{F}(\widetilde{w}_1)), u'(x, T; \widetilde{w}_{1}, \mathfrak{F}(\widetilde{w}_{1}))\right)\in B_{L^2(\Omega_T)}
(u^0,\rho_0) \times B_{H^{-1}(\Omega_T)}(u^1,\rho_1),
\end{equation}
where $\displaystyle B_X(C, r)$ denotes the ball in $X$ with center $C$ and radius $r$.

To explain this optimal problem, we are going to consider the following
sub-problems:

\textbf{$\bullet$ Problem 1} Fixed any leader control $\widetilde{w}_1$,  find the
follower control $\displaystyle \widetilde{w}_2=\mathfrak{F}(\widetilde{w}_1)$ (depending
on $\displaystyle \widetilde{w}_1$) and the associated state $u$, solution of
\eqref{eq1.3.1} satisfying the condition \eqref{son} (Nash equilibrium)  related to
$\widetilde{J}_2$, defined in \eqref{sfn}.

\textbf{$\bullet$ Problem 2} Assuming that the existence of the Nash equilibrium
$\displaystyle \widetilde{w}_2$ was proved, then when $\widetilde{w}_{1}$ varies in
$L^2(\widehat{\Sigma}_1)$, prove that the solutions $\displaystyle \left(u(x, t;
\widetilde{w}_1, \widetilde{w}_2), u'(x, t; \widetilde{w}_1, \widetilde{w}_2)\right)$ of
the state equation \eqref{eq1.3.1}, evaluated at $t = T$, that is, $\displaystyle \left(u(x, T;
\widetilde{w}_1, \widetilde{w}_2), u'(x, T; \widetilde{w}_1, \widetilde{w}_2)\right)$,
generate a dense subset of \break $L^2(\Omega_T) \times H^{-1}(\Omega_T)$.

\begin{remark}\label{r1} By the linearity of system \eqref{eq1.3.1F}, without loss of generality we may assume that $u_0=0=u_1$.\\
\end{remark}

\subsection{Reduction to controllability problem in a cylindrical domain}

 When $0 < k < 1,$ in order to prove the main result of this paper (Theorem \ref{AC}), we first transform \eqref{eq1.3} into a wave equation with variable coefficients in a cylindrical domain.

For this, we divide $\displaystyle {\Sigma}_0^*$ into two
parts
\begin{equation}\label{decomp0.1}
\displaystyle {\Sigma}_0^*=\displaystyle {\Sigma}_1 \cup \displaystyle {\Sigma}_2,
\end{equation}
and consider
\begin{equation} \label{decomp.2}
\displaystyle w = \{w_1,w_2\}, \;\; {w}_i=\mbox{control function in } \; L^2({\Sigma}_i), \;i=1,2.
\end{equation}
We can also write
\begin{equation} \label{decomp 2.A}
\displaystyle w= w _1+ w_2, \; \mbox{ with } \; \displaystyle {\Sigma}_0^*=\displaystyle {\Sigma}_1 = \displaystyle {\Sigma}_2.
\end{equation}

Note that when $\displaystyle (x,t)$ varies in $\displaystyle \widehat{Q}$ the point $\displaystyle (y,t)$,
with $\displaystyle y=\frac{x}{\alpha_k(t)}$, varies in $\displaystyle Q=\Om \times (0,T),$ where
$\Om = (0,1)$. Then, the application $$\displaystyle \zeta:\widehat{Q} \to Q,
\;\;\;\zeta(x,t)=(y,t)$$ is of class $\displaystyle C^2$ and the inverse $\displaystyle \zeta^{-1}$ is
also of class $\displaystyle C^2.$ Observe that in \eqref{decomp 2.A} $\Sigma_i = \zeta(\widehat{\Sigma}_i) (i = 1,2).$ Therefore the change of variables $\displaystyle u(x,t)=v(y,t)$,
transforms the initial-boundary value problem \eqref{eq1.3.1} into the equivalent
system
\begin{equation} \label{eq1.14}
\left|
\begin{array}{l}
\displaystyle v'' + Lv = 0  \ \ \mbox{ in } \ \ Q,\\ [13pt]
\displaystyle v(y,t) = \left\{
\begin{array}{l}
w_1 \ \ \mbox{on} \ \ \Sigma_1,\\
w_2 \ \ \mbox{on} \ \ \Sigma_2,\\
0 \ \ \mbox{on} \ \ \Sigma\backslash\Sigma_0^*,
\end{array}
\right. \\[5pt]
\displaystyle v(y,0) = v_0(y), \; v'(y,0) = v_1(y), \;\; y \in \Om,
\end{array}
\right.
\end{equation}
where
\begin{equation*}
\left|
\begin{array}{l}
\displaystyle Lv = - \Big [\frac{\beta_k(y,t)}{\alpha_k(t)}v_{y}\Big]_{y} + \frac{\gamma_k(y)}{\alpha_k(t)}v'_{y}\;\;, \\[10pt]
\displaystyle \beta_{k}(y,t) = \frac{1  - k^2y^2}{\alpha_{k}(t)},\\[10pt]
\displaystyle \gamma_{k}(y) = -2ky,\\[10pt]
\displaystyle v_{0}(y) = u_{0}(x), \\[10pt]
\displaystyle v_{1}(y) = u_{1}(x) + kyu_{x}(0), \\[10pt]
\displaystyle \Sigma= {\Sigma}_0 \cup \Sigma^*_0,\\[10pt]\displaystyle
\displaystyle {\Sigma}_0 = \{(0,t) : 0 < t < T \},\\[10pt]
\displaystyle {\Sigma}^*_0 = \{(1,t) : 0 < t < T \},\\[10pt]
\displaystyle \displaystyle {\Sigma}_0^* = \Sigma_1  \cup \Sigma_2.
\end{array}
\right.
\end{equation*}

We consider the coefficients of the operator $L$  satisfying  the following
conditions:
\begin{itemize}
\item[(H1)] $\beta_{k}(y,t)  \in C^{1}(\overline{Q});$
\item[(H2)] $\gamma_{k}(y)  \in  W^{1,\infty}(\Om).$
\end{itemize}


\begin{remark}\label{rsol} By a similar method used  in \textup{\cite{Mi}}, we get that for any given $(v_0, v_1) \in L^2(\Om) \times H^{-1}(\Om)$ and  $\displaystyle w_i \in
L^2(\widehat{\Sigma}_i),$ \eqref{eq1.14} admits a unique solution in the sense of transposition,  with  $\displaystyle v \in C\big( [0,T]
;L^2(\Om)) \cap C^{1}\big( [0,T] ; H^{-1}(\Om))$.

Using the diffeomorphism   $\displaystyle \zeta^{-1}(y,t)=(x,t)$, from $Q$ onto
$\widehat{Q}$, we obtain a unique global weak solution $\displaystyle u$ to the
problem \eqref{eq1.3.1} with the regularity \linebreak $\displaystyle u \in C\big( [0,T] ; L^2(\Om_t))
\cap C^{1}\big( [0,T] ; H^{-1}(\Om_t))$.
\end{remark}

$\bullet$ \textbf{Cost functionals in the cylinder $Q$.}  From the diffeomorphism
$\displaystyle \zeta,$ which transforms $\widehat{Q}$ into $Q$, we transform the cost
functionals $\displaystyle \widetilde{J}_2, \displaystyle \widetilde{J}$ into the cost functionals
$J_2, J$ defined by
\begin{equation} \label{eq3.9}
\displaystyle J_2(w_1,w_2) = \displaystyle\frac{1}{2}\,\int_{0}^{T}\int_{\Omega}\alpha_{k}(t)[v(w_1,w_2) - v_2(y,t)]^2dy\,dt + \frac{\sigma}{2}\,\int_{\Sigma_2}w_{2}^{2}\,d\Sigma
\end{equation}
and
\begin{equation} \label{eq3.7}
\displaystyle J(w_1) = \frac{1}{2}\,\int_{\Sigma_1}w_1^{2}\,d\Sigma,
\end{equation}
where $\sigma > 0$ is a constant and $v_2(y,t)$ is a given function in $L^2(\Om
\times (0,T)).$

Associated with the functionals $\displaystyle J_2$ and $\displaystyle J$ defined above, we will
consider the following  sub-problems:

\textbf{$\bullet$ Problem 3} Fixed any leader control ${w}_1$,  find the  follower
control $\displaystyle {w}_2$ (depending on $\displaystyle {w}_1$) and the associated state $v$
solution of \eqref{eq1.14} satisfying (Nash equilibrium)
\begin{equation}\label{soncil}
\displaystyle {J}_2({w}_1, {w}_2)= \inf_{\widehat{w}_2 \in L^2({\Sigma}_2)} {J}_2({w}_1, \widehat{w}_2),
\end{equation}
related to ${J}_2$ defined in \eqref{eq3.9}.

\textbf{$\bullet$ Problem 4} Assuming that the existence of the Nash equilibrium
$\displaystyle {w}_2$ was proved, then when ${w_1}$ varies in $L^2({\Sigma}_1)$, prove
that the solutions $\displaystyle \left(v(y, t; {w_1}, {w}_2), v'(y, t; {w_1}, {w}_2)\right)$ of
the state equation \eqref{eq1.14}, evaluated at $t = T$, that is, $\displaystyle \left(v(y, T;
{w_1}, {w}_2), v'(y, T; {w_1}, {w}_2)\right)$, generate a dense subset of $\displaystyle
L^2(\Om) \times H^{-1}(\Om)$.

The main result of this paper is the following:

\begin{theorem}\label{AC} Assume that \eqref{hT} and \eqref{hT10} hold. Let us consider $\displaystyle w_1 \in L^2(\Sigma_1)$ and $\displaystyle w_2$ a Nash equilibrium in the sense \eqref{soncil}.
Then $\displaystyle \left(v(T), v'(T)\right)=\left(v(., T, {w_1}, w_2), v'(., T, {w_1},
w_2)\right)$, where $\displaystyle v$ solves the system \eqref{eq1.14}, generates a dense
subset of $\displaystyle L^2(\Om)\times H^{-1}(\Om)$.
\end{theorem}

The contents of this paper are organized as follows. Section \ref{sec3} is devoted to establish the existence and uniqueness of Nash equilibrium.  In Section \ref{sec4}, we investigate the approximate controllability proving the density Theorem \ref{AC}. In Section \ref{sec5}, we present the optimality system for the leader control. Finally, we present some additional comments and questions in Section \ref{sec6}.

\section{Nash equilibrium}\label{sec3}
\setcounter{equation}{0}
In this section, fixed any leader control $\displaystyle \ w_1 \in L^2(\Sigma_1)$ we
determine the existence and uniqueness of solutions to the problem
\begin{equation} \label{eq3.10}
\begin{array}{l}
\displaystyle\inf_{w_2 \in L^2(\Sigma_2)}J_2(w_1,w_2),
\end{array}
\end{equation}
and a characterization of this solution in terms of an adjoint system.

In fact, this is a classical type problem in the control of distributed systems (cf.
J.-L. Lions \cite{L3}). It admits a unique solution
\begin{equation} \label{eq3.11}
\displaystyle w_2 = \mathfrak{F}(w_1).
\end{equation}

Indeed, for solution of the problem \eqref{eq3.10}, we will minimize the functional $J_2.$ For this, we consider $\disp \mathcal{U}_{ad} = \{ (v,w_2) \in L^2(Q) \times L^2(\Sigma_2) :  v \ \  \mbox{solution of } \eqref{eq1.14}\}$ and
$\disp J_2 : \mathcal{U}_{ad} \longrightarrow \re$ defined by \eqref{eq3.9}. We write $\disp v=v(w_1,w_2)$ and prove the following items:

\begin{description}

\item[(a)] $\mathcal{U}_{ad}$ is non-empty, closed convex subset of $L^2(Q) \times L^2(\Sigma_2).$

\item[(b)] $J_2$ is weakly coercive.

In fact, using the triangle inequality, we have
\begin{equation*}
\disp ||v-v_2||_{L^2(Q)} \geq \big| ||v||_{L^2(Q)} - ||v_2||_{L^2(Q)}\big|
\end{equation*}
and as $v_2$ is fixed, it follows that
\begin{equation*}
\disp \lim_{{||v||_{L^2(Q)} \rightarrow \infty} \atop {||w_2||_{L^2(\Sigma_2)}\rightarrow \infty}} J_2(v,w_2) = \frac{1}{2}\big|\big|{(\alpha_{k}(t))^\frac{1}{2}}(v-v_2)\big|\big|_{L^2(Q)}^{2} + \frac{\sigma}{2}\,||w_2||_{L^2(\Sigma_2)}^{2} \rightarrow \infty.
\end{equation*}

\item[(c)] $J_2$ is weakly sequentially lower semicontinuous.

Indeed, as $\mathcal{U}_{ad}$ is a closed subset of Hilbert space  $L^2(Q) \times L^2(\Sigma_2)$ then $(v,w_2) \in \mathcal{U}_{ad}.$ We consider two sequences $(v^n) \subset L^2(Q)$ and  $(w_{2}^{n}) \subset  L^2(\Sigma_2)$ such that
\begin{equation*}
\begin{array}{l}
\disp v^n \rightharpoonup v \ \ \ \ \ \mbox{in}\ \ \ \  L^2(Q) \\
\disp w_{2}^{n} \rightharpoonup w_2 \ \ \ \mbox{in}\ \ \  L^2(\Sigma_2).
\end{array}
\end{equation*}
Therefore,
\begin{equation*}
\disp \lim_{n \rightarrow \infty}\inf \frac{1}{2}\big|\big|{(\alpha_{k}(t))^\frac{1}{2}}(v^n-v_2)\big|\big|_{L^2(Q)}  \geq \frac{1}{2} \big|\big|{(\alpha_{k}(t))^\frac{1}{2}}(v-v_2)\big|\big|_{L^2(Q)}
\end{equation*}
and
\begin{equation*}
\disp \lim_{n \rightarrow \infty}\inf ||w_{2}^{n}||_{L^2(\Sigma_2)} \geq ||w_{2}||_{L^2(\Sigma_2)}.
\end{equation*}
Now,
\begin{align*}
\lim_{n \rightarrow \infty}\inf J_2(v^n,w_{2}^{n}) &= \lim_{n \rightarrow \infty}\inf \left\{\frac{1}{2}\big|\big|{(\alpha_{k}(t))^\frac{1}{2}}(v^n-v_2)\big|\big|_{L^2(Q)}^{2}
 + \frac{\sigma}{2}\,||w^n_2||_{L^2(\Sigma_2)}^{2}\right\} \\
&\geq \lim_{n \rightarrow \infty}\inf\left\{\frac{1}{2}\big|\big|{(\alpha_{k}(t))^\frac{1}{2}}(v^n-v_2)\big|\big|_{L^2(Q)}^{2} \right\}\\
 & + \lim_{n \rightarrow \infty}\inf \left\{\frac{\sigma}{2}\,||w^n_2||_{L^2(\Sigma_2)}^{2}\right\}\\
&\geq \frac{1}{2}\big|\big|{(\alpha_{k}(t))^\frac{1}{2}}(v-v_2)\big|\big|_{L^2(Q)}^{2} + \frac{\sigma}{2}||w_2||_{L^2(\Sigma_2)}^{2},\\
\end{align*}
that is,
\begin{equation*}
\lim_{n \rightarrow \infty}\inf J_2(v^n,w_{2}^{n}) \geq J_2(v,w_2).
\end{equation*}

\item[(d)] $J_2$ is strictly convex.

In fact, we consider $\lambda \in (0,1)$ and $(v,w_2)$, $(\widetilde{v}, \widetilde{w}_2) \in  \mathcal{U}_{ad} $ with $(v,w_2) \neq (\widetilde{v}, \widetilde{w}_2)$.
Writing $ v_2 = \lambda\,v_2 + (1 - \lambda)v_2,$ we have
\begin{equation}\label{eq3.72}
\begin{array}{l}
\disp J_2[\lambda(v,w_2) + (1-\lambda)(\widetilde{v}, \widetilde{w}_2)] \\
\disp = \frac{1}{2}\,\int_{0}^{T}\,\int_{\Omega}{\alpha_{k}(t)}[\lambda\,v + (1-\lambda)\widetilde{v} - v_2]^2dy\,dt \\+ \disp \frac{\sigma}{2}\,\int_{\Sigma_2}[\lambda\,w_2 + (1-\lambda)\widetilde{w}_2]^2d\Sigma  \\
\disp = \frac{1}{2}\,\int_{0}^{T}\,\int_{\Omega}{\alpha_{k}(t)}[\lambda\,v + (1-\lambda)\widetilde{v} - \lambda v_2 - (1-\lambda)v_2]^2dy\,dt\\[10pt] + \disp \frac{\sigma}{2}\,\int_{\Sigma_2}[\lambda\,w_2 + (1-\lambda)\widetilde{w}_2]^2d\Sigma \\
\disp = \frac{1}{2}\,\int_{0}^{T}\,\int_{\Omega}{\alpha_{k}(t)}[\lambda(v-v_2) + (1-\lambda)(\widetilde{v}-v_2)]^2dy\,dt\\[10pt] \disp + \frac{\sigma}{2}\,\int_{\Sigma_2}[\lambda\,w_2 + (1-\lambda)\widetilde{w}_2]^2d\Sigma.
\end{array}
\end{equation}

Expanding the expression after the last equal sign in (\ref{eq3.72}), we obtain
\begin{equation} \label{eq3.73}
\begin{array}{l}
\disp \frac{\lambda^2}{2}\,\int_{0}^{T}\int_{\Omega}{\alpha_{k}(t)}(v-v_2)^2dy\,dt\\[10pt] \disp  + \lambda(1-\lambda)\int_{0}^{T}\int_{\Omega}{\alpha_{k}(t)}\underbrace{(v-v_2)(\widetilde{v}-v_2)}_{(\ast)} dy\,dt \\[10pt]
\disp + \frac{(1-\lambda)^2}{2}\int_{0}^{T}\int_{\Omega}{\alpha_{k}(t)}(\widetilde{v}-v_2)^2dy\,dt + \frac{\sigma\,\lambda^2}{2}\int_{\Sigma_2}w_{2}^{2}\, d\Sigma \\[10pt] \disp + \sigma\,\lambda(1-\lambda)\int_{\Sigma_2}\underbrace{w_2\widetilde{w}_2}_{(\ast\ast)}\,d\Sigma +\frac{\sigma(1-\lambda)^2}{2}\int_{\Sigma_2}\widetilde{w}_{2}^{2}\,d\Sigma.
\end{array}
\end{equation}

After applying Young's inequality to $(\ast)$ and noting that $\lambda(1 - \lambda) > 0,$ it follows that
\begin{equation} \label{eq3.74}
\begin{array}{l}
\disp \lambda(1-\lambda)\int_{0}^{T}\int_{\Omega}{\alpha_{k}(t)}(v-v_2)(\widetilde{v}-v_2)dydt \\[10pt]
\disp < \frac{\lambda(1-\lambda)}{2}\int_{0}^{T}\int_{\Omega}{\alpha_{k}(t)}(v-v_2)^2dydt \\ \disp + \frac{\lambda(1-\lambda)}{2}\int_{0}^{T}\int_{\Omega}{\alpha_{k}(t)}(\widetilde{v}-v_2)^2dydt.
\end{array}
\end{equation}

Similarly, we estimate  $(\ast\ast)$ as
\begin{equation} \label{eq3.75}
\begin{array}{l}
\disp \sigma\,\lambda(1-\lambda)\int_{\Sigma_2}w_2\widetilde{w}_2\,d\Sigma\, < \, \frac{\sigma\,\lambda(1-\lambda)}{2}\int_{\Sigma_2}w_{2}^{2}\,d\Sigma \\ \disp + \frac{\sigma\,\lambda(1-\lambda)}{2}\int_{\Sigma_2}\widetilde{w}_{2}^{2}\,d\Sigma.
\end{array}
\end{equation}

Substituting  (\ref{eq3.74}) and (\ref{eq3.75}) into the right-hand side of (\ref{eq3.72}), we get
\begin{equation*}
\begin{array}{l}
\disp J_2[\lambda(v,w_2) + (1-\lambda)(\widetilde{v},\widetilde{w}_2)] \\[10pt]
\disp <  \frac{\lambda^2}{2}\int_{0}^{T}\int_{\Omega}{\alpha_{k}(t)}(v-v_2)^2 dy\,dt \\[10pt]
\disp + \frac{\lambda(1-\lambda)}{2}\int_{0}^{T}\int_{\Omega}{\alpha_{k}(t)}(v-v_2)^2 dy\,dt \\[10pt]
\disp + \frac{\lambda(1-\lambda)}{2}\int_{0}^{T}\int_{\Omega}{\alpha_{k}(t)}(\widetilde{v}-v_2)^2 dy\,dt  \\[10pt]
\disp + \frac{(1-\lambda)^2}{2}\int_{0}^{T}\int_{\Omega}{\alpha_{k}(t)}(\widetilde{v}-v_2)^2 dy\,dt
+ \frac{\sigma\lambda^2}{2}\int_{\Sigma_2}w_{2}^{2}\,d\Sigma \\[10pt]
\disp + \frac{\sigma\lambda(1-\lambda)}{2}\int_{\Sigma_2}w_{2}^{2}\,d\Sigma
\disp + \frac{\sigma\lambda(1-\lambda)}{2}\int_{\Sigma_2}\widetilde{w}_{2}^{2}\,d\Sigma \\[10pt]
+ \disp \frac{\sigma(1-\lambda)^2}{2}\int_{\Sigma_2}\widetilde{w}_{2}^{2}\,d\Sigma = \frac{\lambda}{2}\int_{0}^{T}\int_{\Omega}{\alpha_{k}(t)}(v-v_2)^2 dy\,dt
+ \frac{\sigma\lambda}{2}\int_{\Sigma_2}w_{2}^{2}\,d\Sigma \\[10pt]
\disp + \frac{(1-\lambda)}{2}\int_{0}^{T}\int_{\Omega}{\alpha_{k}(t)}(\widetilde{v}-v_2)^2 dy\,dt +
\frac{\sigma(1-\lambda)}{2}\int_{\Sigma_2}\widetilde{w}_{2}^{2}\,d\Sigma  \\[10pt]

\disp = \lambda\left[\frac{1}{2}\int_{0}^{T}\int_{\Omega}{\alpha_{k}(t)}(v-v_2)^2 dy\,dt + \frac{\sigma}{2}\int_{\Sigma_2}w_{2}^{2}\,d\Sigma\right]  \\[10pt]
\disp + (1-\lambda)\left[ \frac{1}{2}\int_{0}^{T}\int_{\Omega}{\alpha_{k}(t)}(\widetilde{v}-v_2)^2 dy\,dt +  \frac{\sigma}{2}\int_{\Sigma_2}\widetilde{w}_{2}^{2}\,d\Sigma\right] \\[10pt]
\disp = \lambda\,J_2(v,w_2) + (1-\lambda)J_2(\widetilde{v},\widetilde{w}_2),
\end{array}
\end{equation*}
that is,
\begin{equation*}
\disp J_2[\lambda(v,w_2) + (1-\lambda)(\widetilde{v}, \widetilde{w}_2)] < \lambda J_2(v,w_2) + (1-\lambda)J_2(\widetilde{v}, \widetilde{w}_2).
\end{equation*}
\end{description}

Now, we will calculate the  Gateaux derivative of the functional \eqref{eq3.9}. For $\theta_1 \in L^2\big(\Omega \times (0,T)\big),$  $\theta_2 \in L^2\big(\Sigma_2\big)$ and $\varepsilon > 0$, we have
\begin{align*}
&\nonumber J_{2}'(v,w_2) = \lim_{\varepsilon \rightarrow 0}\frac{1}{\varepsilon}\bigg\{J_{2}(v + \varepsilon\theta_1, w_2 + \varepsilon\theta_2) - J_{2}(v, w_2)\bigg\} \\
\nonumber& = \lim_{\varepsilon \rightarrow 0}\frac{1}{\varepsilon}\Bigg\{ \frac{1}{2}\int_{0}^{T}\int_{\Omega}{\alpha_{k}(t)}(v+\varepsilon\theta_1-v_2)^2 dy\,dt  + \frac{\sigma}{2}\int_{\Sigma_2}(w_{2}   + \varepsilon\theta_2)^2d\Sigma \\ &-
\frac{1}{2}\int_{0}^{T}\int_{\Omega}{\alpha_{k}(t)}(v-v_2)^2 dy\,dt - \frac{\sigma}{2}\int_{\Sigma_2}w_{2}^{2}\,d\Sigma \Bigg\} \\
\nonumber&= \lim_{\varepsilon \rightarrow 0}\frac{1}{\varepsilon}\Bigg\{\frac{1}{2}\int_{0}^{T}\int_{\Omega}{\alpha_{k}(t)}\Big[(v-v_2)^2 + 2\varepsilon\theta_1(v-v_2) + \varepsilon^2\theta_{1}^{2}\Big] dy\,dt\\ &+ \frac{\sigma}{2}\int_{\Sigma_2}(w_{2} + \varepsilon\theta_2)^2d\Sigma
\nonumber - \frac{1}{2}\int_{0}^{T}\int_{\Omega}{\alpha_{k}(t)}(v-v_2)^2 dy\,dt - \frac{\sigma}{2}\int_{\Sigma_2}w_{2}^{2}\,d\Sigma \Bigg\} \\&= \lim_{\varepsilon \rightarrow 0}\frac{1}{\varepsilon}\Bigg\{\frac{1}{2}\int_{0}^{T}\int_{\Omega}{\alpha_{k}(t)}(v-v_2)^2dy\,dt
\nonumber+ \varepsilon\int_{0}^{T}\int_{\Omega}{\alpha_{k}(t)}(v-v_2)\theta_1dy\,dt\\ &+ \frac{1}{2}\,\varepsilon^2\int_{0}^{T}\int_{\Omega}{\alpha_{k}(t)}\theta_{1}^{2}dy\,dt + \frac{\sigma}{2}\int_{\Sigma_2}w_{2}^{2}d\Sigma + \varepsilon\sigma\int_{\Sigma_2}w_2\theta_2d\Sigma
+ \frac{\sigma}{2}\,\varepsilon^2\int_{\Sigma_2}\theta^2d\Sigma \\&- \frac{1}{2}\int_{0}^{T}\int_{\Omega}{\alpha_{k}(t)}(v-v_2)^2 dy\,dt - \frac{\sigma}{2}\int_{\Sigma_2}w_{2}^{2}\,d\Sigma \Bigg\} \\
& = \int_{0}^{T}\int_{\Omega}{\alpha_{k}(t)}(v-v_2)\theta_1dy\,dt + \sigma\int_{\Sigma_2}w_2\theta_2d\Sigma.
\end{align*}

Therefore, the Euler - Lagrange equation for problem  \eqref{eq3.10} is given by
\begin{equation} \label{eq3.21}
\int_{0}^{T}\int_{\Omega}{\alpha_{k}(t)}(v-v_2)\widehat{v}dy\,dt + \sigma\int_{\Sigma_2}w_2\widehat{w}_2d\Sigma = 0, \;\;\forall\, \widehat{w}_2 \in L^2(\Sigma_2),
\end{equation}
where $\widehat{v}$ is solution of the following system
\begin{equation} \label{eq3.22}
\left|
\begin{array}{l}
\displaystyle \widehat v'' +L\widehat{v} = 0 \ \ \mbox{ in } \ \ Q,\\ [7pt]
\displaystyle \widehat{v} = \left\{
\begin{array}{l}
0 \ \ \mbox{ on } \ \ \Sigma_1, \\
\widehat{w}_2 \ \ \mbox{ on } \ \ \Sigma_2, \\
0 \ \ \mbox{ on } \ \ \Sigma \backslash \left(\Sigma_1 \cup \Sigma_2\right),
\end{array}
\right. \\[13pt]
\displaystyle \widehat{v}(y,0) = 0, \; \widehat{v'}(y,0) = 0, \;\; y \in \Om.
\end{array}
\right.
\end{equation}
In order to express \eqref{eq3.21} in a convenient form, we introduce the adjoint
state defined by
\begin{equation}\label{sac}
\left|
\begin{array}{l}
p'' + L^{\ast}\,p = \alpha_{k}(t)\left(v - v_2\right) \ \mbox{ in } \ \ Q, \\[5pt]\displaystyle
p(T) = p'(T) = 0, \;\; y \in \Om, \\[5pt]\displaystyle
p = 0 \ \mbox{ on } \ \Sigma,
\end{array}
\right.
\end{equation}
where $L^{\ast}$ is the formal adjoint of the  operator $\displaystyle L$.

Multiplying \eqref{sac} by $\widehat{v}$ and integrating by parts, we find
\begin{equation} \label{eq3.33}
\int_{0}^{T}\int_\Om \alpha_{k}(t)(v - v_2)\widehat{v}\,dy\,dt + \int_{\Sigma_2} \frac{1}{\alpha_k^2(t)}\,p_y\,\widehat{w}_2\,d\Sigma = 0,
\end{equation}
so that \eqref{eq3.21} becomes
\begin{equation}\label{ci}
\displaystyle p_y= \sigma \alpha_k^2(t)\,w_2 \ \ \mbox{ on } \ \ \Sigma_2.
\end{equation}
We summarize these results in the following theorem.
\begin{theorem}\label{teN} For each $\displaystyle w_1 \in L^2(\Sigma_1)$ there exists a unique Nash equilibrium $\displaystyle w_2$ in the sense of \eqref{soncil}. Moreover, the follower $\displaystyle w_2$ is given by
\begin{equation}\label{cseg}
\displaystyle \displaystyle w_2 = \mathfrak{F}(w_1)=\frac{1}{\sigma \alpha_{k}^2(t)}\,\;p_y\;\;\mbox{ on }\;\;\Sigma_2,
\end{equation}
where $\displaystyle \{ v,p \}$ is the unique solution of (the optimality system)
\begin{equation} \label{eq3.37}
\left|
\begin{array}{l}
\displaystyle v'' + Lv = 0 \ \mbox{ in } \;\; Q, \\
\displaystyle p'' + L^{\ast}\,p = \alpha_{k}(t)\left(v - v_2\right) \ \mbox{ in } \;\; Q,\\[5pt]
\displaystyle v = \left\{
\begin{array}{l}
w_1 \ \mbox{ on } \ \Sigma_1,\\[5pt]
\displaystyle \frac{1}{\sigma \alpha_{k}^2(t)}\;\,p_y \ \mbox{ on } \ \Sigma_2,\\[5pt]
0 \ \mbox{ on } \ \Sigma \backslash \Sigma_0^*,
\end{array}
\right.\\[7pt]
p = 0 \ \mbox{ on } \ \Sigma, \\[5pt]
v(0) = v'(0) = 0, \\[5pt]
p(T) = p'(T) = 0, \;\; y \in \Om.
\end{array}
\right.
\end{equation}
Of course, $\displaystyle \{ v,p \}$ depends on $w_1$:
\begin{equation}\label{cdep}
\displaystyle \{ v,p \} = \{ v(w_1),p(w_1) \}.
\end{equation}
\end{theorem}
\section{Approximate controllability}\label{sec4}
\setcounter{equation}{0}
Since we have proved the existence, uniqueness and characterization of the
follower $\displaystyle w_2$, the leader $\displaystyle w_1$  now wants  that solutions $v$
and $v'$, evaluated at time $t=T$, to  be as close as possible to $\displaystyle (v^0,
v^1)$. This will be possible if the system \eqref{eq3.37}  is approximately
controllable. We are looking for
\begin{equation} \label{inf1}
\begin{array}{l}
\displaystyle\inf\, \frac{1}{2\,}\,\int_{\Sigma_1} w_{1}^{2}\,d\Sigma,
\end{array}
\end{equation}
where $\displaystyle w_1$ is subject to
\begin{equation} \label{subj1}
\begin{array}{l}
\displaystyle \left(v(T;{w_1}), v'(T; {w_1})\right)\in B_{L^2(\Om)}(v^0,\rho_0) \times B_{H^{-1}(\Om)}(v^1,\rho_1),
\end{array}
\end{equation}
assuming that  $w_1$ exists, $\rho_0$ and  $\rho_1$ being positive numbers arbitrarily
small  and $\{v^0, v^1\} \in L^2(\Om) \times H^{-1}(\Om)$.

As in \cite{Cui1}, we assume that
\begin{equation}\label{hT}
\disp T >  \frac{e^{\frac{2k(1 + k)}{1 - k}} - 1}{k}
\end{equation}
and
\begin{equation}\label{hT10}
0 < k <  1.
\end{equation}

Now as in the case \eqref{decomp 2.A} and using  Holmgren's Uniqueness
Theorem (cf. \cite{LH}; and see also \cite {Cui1} for additional discussions), we conclude this section with the proof of Theorem \ref{AC}.

\noindent{\bf{Proof of Theorem \ref{AC}}}

We decompose the solution $\displaystyle (v,p)$ of \eqref{eq3.37}
setting
\begin{equation} \label{eq3.39}
\left|
\begin{array}{l}
v = v_0 + g,\\
p = p_0 + q,
\end{array}
\right.
\end{equation}
where $v_0$, $p_0$ is given by
\begin{equation} \label{eq3.40}
\left|
\begin{array}{l}
\displaystyle v_0'' + L\,v_0 = 0 \ \mbox{ in } \ Q,\\[5pt]\displaystyle
v_0 = \left\{
\begin{array}{l}
0 \ \mbox{ on } \ \Sigma_1,\\[5pt]
\displaystyle \frac{1}{\sigma \alpha_{k}^2(t)}\,({p_0})_{y} \ \mbox{ on } \ \Sigma_2,\\[5pt]\displaystyle
0 \ \mbox{ on } \ \Sigma \backslash \Sigma_0^*,
\end{array}
\right.\\[5pt]\displaystyle
v_0(0) = v_0'(0) = 0, \;\; y\in \Om,
\end{array}
\right.
\end{equation}
\begin{equation} \label{eq3.41}
\left|
\begin{array}{l}
\displaystyle p_0'' + L^{\ast}p_0 = \alpha_{k}(t) \left(v_0 - v_2\right) \ \mbox{ in } \ Q, \\[5pt]\displaystyle
p_0 = 0 \ \mbox{ on } \ \Sigma,\\[5pt]\displaystyle
p_0(T) = p_0'(T) = 0, \;\; y \in \Om,
\end{array}
\right.
\end{equation}
and $\displaystyle \{g,q\}$ is given by
\begin{equation} \label{eq3.42}
\left|
\begin{array}{l}
g'' + L\,g = 0 \ \mbox{ in } \ Q,\\[5pt]\displaystyle
g = \left\{
\begin{array}{l}
w_1 \ \mbox{ on } \ \Sigma_1,\\[5pt]
\displaystyle \frac{1}{\sigma \alpha_{k}^2(t)}\,q_{y} \ \mbox{ on } \ \Sigma_2,\\[5pt]\displaystyle
0 \ \mbox{on} \ \Sigma \backslash \Sigma_0^*,
\end{array}
\right.\\[5pt]\displaystyle
g(0) = g'(0) = 0, \;\;y \in \Om,
\end{array}
\right.
\end{equation}
\begin{equation} \label{eq3.43}
\left|
\begin{array}{l}
\displaystyle q'' + L^{\ast}q = \alpha_{k}(t) g \ \mbox{ in } \ Q, \\[5pt]\displaystyle
q = 0 \ \mbox{ on } \ \Sigma,\\[5pt]\displaystyle
q(T) = q'(T) = 0, \;\; y \in \Om.
\end{array}
\right.
\end{equation}
We next set
\begin{equation} \label{eq3.44}
\begin{array}{ccll}
A \ : & \! L^2(\Sigma_1) & \! \longrightarrow & \! H^{-1}(\Om) \times L^2(\Om) \\
& \! w_1 & \! \longmapsto & \! A\,w_1 = \big\{ g'(T;w_1) + \delta g(T;w_1),\; -g(T;w_1) \big\},
\end{array}
\end{equation}
which defines
$$A \in \mathcal{L}\left( L^2(\Sigma_1); \;H^{-1}(\Om) \times L^2(\Om)\right),$$
where  $\delta$ is a positive constant.

Using \eqref{eq3.39} and \eqref{eq3.44}, we can rewrite  \eqref{subj1} as
\begin{equation} \label{subj2}
\displaystyle Aw_1\in \{ -v_0(T)+\delta g(T)+B_{H^{-1}(\Om)}(v^1,\rho_1),\;-v_0(T)+B_{L^2(\Om)}(v^0,\rho_0)\}.
\end{equation}
We will show that  $Aw_1$ generates a dense subspace of $H^{-1}(\Om) \times
L^2(\Om)$. For this, let  $\{ f^0,f^1 \} \in H_{0}^{1}(\Om) \times L^2(\Om)$ and
consider the following systems (``adjoint states"):
\begin{equation} \label{eq3.45}
\left|
\begin{array}{l}
\varphi'' + L^{\ast}\,\varphi =\displaystyle \alpha_{k}(t) \;\psi \ \mbox{ in } \ Q, \\[5pt]\displaystyle
\varphi = 0 \ \mbox{ on } \ \Sigma,\\[5pt]\displaystyle
\varphi(T) = f^0, \ \varphi'(T) = f^1, \;\; y \in \Om,
\end{array}
\right.
\end{equation}

\begin{equation} \label{eq3.46}
\left|
\begin{array}{l}
\psi'' + L\,\psi = 0 \ \mbox{ in } \ Q,\\[5pt]\displaystyle
\psi = \left\{
\begin{array}{l}
0 \ \mbox{ on } \ \Sigma_1,\\[5pt]\displaystyle
\displaystyle \frac{1}{\sigma \alpha_{k}^2(t)}\,\varphi_{y} \ \mbox{ on } \ \Sigma_2,\\[5pt]\displaystyle
0 \ \mbox{ on } \ \Sigma \backslash \Sigma_0^*,
\end{array}
\right.\\[5pt]\displaystyle
\psi(0) = \psi'(0) = 0, \;\; y \in \Om.
\end{array}
\right.
\end{equation}
Multiplying \eqref{eq3.46}$_1$ by $q$, \eqref{eq3.45}$_1$ by $g$, where $q$, $g$
solve (\ref{eq3.43}) and \eqref{eq3.42}, respectively, and integrating in $Q$ we
obtain
\begin{equation} \label{eq3.47}
\int_{0}^{T}\int_{\Om} \alpha_{k}(t)g\,\psi\,dy\,dt =- \frac{1}{\sigma}\int_{\Sigma_2} \frac{1}{\alpha_{k}^4(t)}\,q_{y}\, \varphi_{y} d\Sigma,
\end{equation}
and
\begin{equation} \label{eq3.49}
\begin{array}{l}
\displaystyle \langle g'(T),f^0 \rangle_{H^{-1}(\Om) \times  H_{0}^{1}(\Om)} + \delta \langle g(T), f^0 \rangle_{L^2(\Om) \times  H_{0}^1(\Om)}  - \big( g(T),f^1 \big)\\ = \displaystyle -\int_{\Sigma_1}\frac{1}{\alpha_{k}^2(t)}\,\varphi_{y}\,w_1\,d\Sigma.
\end{array}
\end{equation}
Considering the left-hand side of the above equation as the inner product of $\displaystyle \{g'(T)+
\delta g(T),-g(T)\}$ with $\{ f^0,f^1 \}$ in $ H^{-1}(\Om)  \times L^2(\Om) $ and $
H_{0}^{1}(\Om) \times L^2(\Om)$, we obtain

\begin{equation*}
\Big\langle \big\langle A\,w_1 , f \big\rangle \Big\rangle = - \int_{\Sigma_1}\frac{1}{\alpha_{k}^2(t)}\,\varphi_{y}\,w_1\,d\Sigma,
\end{equation*}
where $\Big\langle \big\langle . , . \big\rangle \Big\rangle$ represent the duality
pairing between $ H^{-1}(\Om) \times L^2(\Om) $ and $ H_{0}^{1}(\Om) \times
L^2(\Om) $. Therefore, if
\begin{equation*}
\begin{array}{l}
 \disp \Big\langle \big\langle A\,w_1 , f \big\rangle \Big\rangle  = 0
\end{array}
\end{equation*}
for all $w_1 \in L^2(\Sigma_1)$, then
\begin{equation} \label{eq3.50}
\displaystyle \varphi_{y}= 0 \ \ \mbox{ on } \ \ \Sigma_1.
\end{equation}
Hence,  in case \eqref{decomp 2.A},
\begin{equation} \label{eq3.51}
\psi = 0 \ \ \mbox{ on } \ \ \Sigma, \;\; \mbox{ so that } \psi\equiv 0.
\end{equation}
Therefore
\begin{equation} \label{eq3.54}
\begin{array}{l}
\displaystyle \varphi'' + L^{\ast}\,\varphi = 0, \;\; \varphi = 0 \mbox{ on }  \Sigma,
\end{array}
\end{equation}
and satisfies \eqref{eq3.50}.  Therefore, according to  Holmgren's Uniqueness
Theorem (cf. \cite{LH}; and see also \cite {Cui1} for additional discussions) and if \eqref{hT} holds, then $\displaystyle \varphi \equiv 0$,
so that $\displaystyle f^0=0, f^1=0$. This
ends the proof. \cqdf

\section{ The optimality system}\label{sec5}
\setcounter{equation}{0}
Thanks to the results obtained in Section \ref{sec3} , we can take, for  each
$\displaystyle w_1$, the Nash equilibrium $\displaystyle w_2$ associated to solution $\displaystyle v$
of \eqref{eq1.14}. We will show the existence of a leader control $\displaystyle w_1$ solution
of the following problem:
\begin{equation} \label{eq3.7cil.1}
\displaystyle \inf J(w_1),
\end{equation}
where $J(w_1)$ is given in \eqref{eq3.7}. For this, we will use  a duality argument due to Fenchel and Rockafellar \cite{R}
(cf. also \cite{Bre, EK}).

The following result holds:
\begin{theorem} \label{teor3.6} Assume the hypotheses $\displaystyle (H1) - (H2)$,
\eqref{decomp 2.A}, \eqref{hT} and \eqref{hT10}  are satisfied. Then for $\{f^0,f^1\}$
in $\displaystyle H_0^1(\Om) \times L^2(\Om)$ we uniquely define $\{\varphi, \psi, v, p \}$
by
\begin{equation} \label{eq3.139}
\left|
\begin{array}{l}
\displaystyle \varphi'' + L^* \varphi = \alpha_{k}(t)\psi \ \ \text{in} \ \ Q, \\[3pt]\displaystyle
\displaystyle \psi'' + L \psi = 0 \ \ \text{in} \ \ Q, \\[3pt]\displaystyle
\displaystyle v'' + Lv = 0 \ \ \text{in} \ \ Q, \\[3pt]\displaystyle
\displaystyle p'' + L^*p  = \alpha_{k}(t)(v - v_2) \ \ \text{in} \ \ Q, \\[3pt]\displaystyle
\displaystyle \varphi = 0 \ \ \ \text{on} \ \ \ \Sigma, \\[3pt]\displaystyle
\displaystyle \psi =
\left\{
\begin{array}{l}
\displaystyle 0 \ \ \text{on} \ \ \Sigma_1, \\[3pt]\displaystyle
\displaystyle \frac{1}{\sigma \alpha_{k}^2(t)}\;\varphi_{y} \ \ \text{on} \ \ \Sigma_2, \\[3pt]\displaystyle
\displaystyle 0 \ \ \ \text{on} \ \ \ \Sigma \backslash \Sigma_0^*,\\[3pt]\displaystyle
\end{array}
\right.\\
\displaystyle v =
\left\{
\begin{array}{l}
\displaystyle - \frac {1}{\alpha_{k}^2(t)}\;\varphi_{y} \ \ \text{on} \ \ \Sigma_1,\\[3pt]\displaystyle
\displaystyle \frac{1}{\sigma \alpha_{k}^2(t)}\;p_{y} \ \ \text{on} \ \ \Sigma_2,\\[3pt]\displaystyle
\displaystyle 0 \ \ \ \text{on} \ \ \ \Sigma \backslash \Sigma_0^*,\\[3pt]\displaystyle
\end{array}
\right.\\
\displaystyle p = 0 \ \ \ \text{on} \ \ \ \Sigma, \\[3pt]\displaystyle
\displaystyle \varphi(.,T) = f^0,\, \varphi'(.,T) = f^1 \ \ \text{in} \ \  \Om, \\[3pt]\displaystyle
\displaystyle v(0) = v'(0) = 0 \ \ \text{in} \ \   \Om, \\[3pt]\displaystyle
\displaystyle p(T) = p'(T) = 0 \ \ \text{in} \ \  \Om.
\end{array}
\right.
\end{equation}
We uniquely define  $\{f^0,f^1\}$ as the solution of the variational inequality
\begin{equation} \label{eq3.140}
\begin{array}{l}
\displaystyle \big\langle v'(T,f) - v^1, \widehat{f}^0 - f^0\big\rangle_{H^{-1}(\Om) \times  H_{0}^{1}(\Om)} - \big(v(T,f) - v^0, \widehat{f}^1 - f^1\big)  \\[10pt]
\displaystyle + \rho_1\big(||\widehat{f}^0|| - ||f^0||\big) + \rho_0\big(|\widehat{f}^1| - |f^1|\big) \geq 0,\,\forall\, \widehat{f} \in H_{0}^{1}(\Om) \times L^2(\Om).
\end{array}
\end{equation}
Then the optimal leader is given by
\begin{equation*}
w_1 = -\frac{1}{\alpha_{k}^2(t)}\;\varphi_{y} \ \ \text{on} \ \ \Sigma_1,
\end{equation*}
where $\varphi$ corresponds to the solution of \eqref{eq3.139}.
\end{theorem}
\textbf{Proof.} Let $A$ be the continuous linear operator defined by \eqref{eq3.44} and introduce the following two convex proper functions:
\begin{equation}\label{eq3.119}
\begin{array}{l}
\displaystyle F_1 : L^2(\Sigma_1) \longrightarrow \re \cup \{\infty\},\\[5pt]
\displaystyle F_1(w_1) = \frac{1}{2} \int_{\Sigma_1}w_{1}^{2}\,d\Sigma
\end{array}
\end{equation}
and
\begin{equation*}
F_2 : H^{-1}(\Omega) \times L^2(\Omega) \longrightarrow \re \cup \{\infty\},
\end{equation*}
given by
\begin{equation} \label{eq3.120}
\begin{array}{ccl}
F_2(Aw_1) = F_2\big(\{g'(T,w_1) + \delta g(T,w_1),-g(T,w_1)\}\big) \\
= \left\{
\begin{array}{l}
0, \text{ if }
\left\{
\begin{array}{l}
g'(T) + \delta g(T) \in v^1 - v_0'(T) + \delta g(T) + \rho_1B_{H^{-1}(\Om)},\\
-g(T) \in -v^0 + v_0(T,w_1) - \rho_0B_{L^2(\Om)},
\end{array}
\right.\\
+ \infty, \text{ otherwise}.
\end{array}
\right.
\end{array}
\end{equation}
With these notations, problems \eqref{inf1}--\eqref{subj1} become equivalent to
\begin{equation} \label{eq3.122}
\begin{array}{l}
\displaystyle \inf_{w_1 \in L^2(\Sigma_1)}\big[F_1(w_1) + F_2(Aw_1)\big]
\end{array}
\end{equation}
provided we prove that the range of $\displaystyle A$ is dense in $\displaystyle H^{-1}(\Om)
\times L^2(\Om)$, under conditions \eqref{hT} and  \eqref{hT10}.

By the Duality Theorem of Fenchel and Rockafellar \cite{R}(see also \cite{Bre}, \cite{EK}), we have
\begin{equation} \label{eq3.124}
\begin{array}{l}
\displaystyle \inf_{w_1 \in L^2(\Sigma_1)}[F_1(w_1) + F_2(Aw_1)]\\[5pt]\displaystyle  = -\inf_{(\widehat{f}^0,\widehat{f}^1) \in H_{0}^{1}(\Om) \times L^2(\Om)} [F_{1}^{*}\big(A^*\{\widehat{f}^0,\widehat{f}^1\}\big) + F_{2}^{*}\{-\widehat{f}^0, -\widehat{f}^1\}],
\end{array}
\end{equation}
where $\displaystyle F_i^*$ is the conjugate function of $\displaystyle F_i  (i=1,2)$ and $\displaystyle A^*$ the adjoint of $\displaystyle A$.

We have
\begin{equation} \label{eq3.121}
\begin{array}{ccccc}
A^* \ : & \! H_{0}^{1}(\Omega) \times L^2(\Omega) & \! \longrightarrow & \! L^2(\Sigma_1) \\
& \! (f^0,f^1) & \! \longmapsto & \! A^*f = & \! -\dfrac{1}{\alpha_{k}^2(t)}\,\varphi_{y},
\end{array}
\end{equation}
where $\varphi$ is given in (\ref{eq3.45}).
We see easily that
\begin{equation} \label{eq3.125}
F_{1}^{*}(w_1) = F_1(w_1)
\end{equation}
and
\begin{equation}\label{eq3.125.2}
\begin{array}{ccl}
\displaystyle  F_{2}^{*}(\{\widehat{f}^0,\widehat{f}^1 \})  &=& \displaystyle \langle v^1 - v_0'(T) + \delta g(T), \widehat{f}^0\rangle_{H^{-1}(\Omega) \times H_{0}^{1}(\Omega)} \\&&+\big( v_0(T) - v^0 ,\widehat{f}^1\big)
\displaystyle + \rho_1||\widehat{f}^0|| + \rho_0|\widehat{f}^1|.
\end{array}
\end{equation}

So, from \eqref{eq3.121}--\eqref{eq3.125.2}, the expression in \eqref {eq3.124}  becomes
\begin{equation} \label{eq3.1240}
\begin{array}{l}
\displaystyle \inf_{w_1 \in L^2(\Sigma_1)}[F_1(w_1) + F_2(Aw_1)] = -\inf_{(\widehat{f}^0,\widehat{f}^1) \in  H_{0}^{1}(\Om) \times L^2(\Om)}  \Theta \big(\{\widehat{f}^0,\widehat{f}^1\}\big),
\end{array}
\end{equation}
where the functional $\Theta : H_{0}^{1}(\Omega) \times L^2(\Omega) \longrightarrow \re$ is defined by
\begin{equation*}
\begin{array}{l}
\disp \Theta \big(\{\widehat{f}^0,\widehat{f}^1\}\big) = \frac{1}{2}\int_{\Sigma_1}\left(\frac{1}{\alpha_{k}^2(t)}\right)^2  \widehat{\varphi}_{y}^2d\;\Sigma  + \big( v^0 - v_0(T) ,\widehat{f}^1\big) \\
 - \langle v^1 - v_0'(T) + \delta g(T), \widehat{f}^0\rangle_{ H^{-1}(\Omega) \times  H_{0}^{1}(\Omega)} + \rho_1||\widehat{f}^0|| + \rho_0|\widehat{f}^1|.
\end{array}
\end{equation*}
This is the dual problem of  \eqref{inf1},  \eqref{subj1}.  Hence, we can use the primal or the dual problem to derive the optimality  system for the leader control. \cqdf

\section{Some additional comments and questions}\label{sec6}

As a future work we are looking for improvements and generalizations of these results to other models. To close this section, we make some comments and briefly discuss some possible extensions of our results and also indicate open issues on the subject.

\textbf{$\bullet$} In the case of $k = 1,$ some results have been obtained in \cite {Cui1}. However, we do not extend the approach developed in this paper to the case $k = 1.$

\textbf{$\bullet$} In the case  $k > 1,$ the moving boundary is a spacelike surface, on which an initial condition rather than a boundary condition needs to be imposed.  It would be quite interesting to study the controllability for \eqref{eq1.3} in this case.  For interested readers on this subject, we cite for instance \cite {Cui1}, \cite {Cui100}, and \cite {SUN}.

\textbf{$\bullet$}  In \cite{AR501}, those authors considered the controllability for semilinear wave equations using a fixed-point formulation. Moreover, they proposed as a future work possible extensions to boundary controls and/or superlinear nonlinearities.

\paragraph{\bf Acknowledgements}

The author wants to express his gratitude to the anonymous reviewers for their questions and commentaries; they were very
helpful in improving this article.

 \paragraph
\end{document}